# A HETEROPOLYMER IN A MEDIUM WITH RANDOM DROPLETS

By Mario V. Wüthrich

*ETH Zürich*

We define a heteropolymer in a medium with random droplets. We prove that for this model we have two regimes: a delocalized one and a localized one. In the localized regime we prove tightness to the droplets, whereas in the delocalized regime we prove diffusive path behavior.

## 1. Introduction and main result.

1.1. *Definition of the model.* We start with the definition of the model. The dimension is denoted by $d+1$, $d \geq 1$. Then our model has three random ingredients: (1) a $(d+1)$-dimensional random medium, described by $\eta$; (2) a $(d+1)$-dimensional directed random walk, described by $S$; (3) a sequence of random variables $\omega$, taking values $\pm 1$ describing the random monomers. So let us formally define all these objects.

1. *Random medium.* Choose $\eta = (\eta_i)_{i \in \mathbb{Z}}$: an i.i.d. sequence of random variables taking the value 1 with probability $p \in [0,1]$ and taking value $-1$ with probability $1-p$. Then the random medium is defined as follows:

$$(1) \qquad D_\eta = \{(k, 0, \ldots, 0) \in \mathbb{Z}^{d+1} \text{ if } \eta_k = 1\}.$$

$D_\eta$ describes a $(d+1)$-dimensional medium with random droplets attached to the first coordinate axis, that is, we have a droplet in $(k, 0, \ldots, 0)$ whenever $\eta_k = 1$. The probability measure of $\eta$ is denoted by $\mathbb{P}_\eta$.

2. *Directed random walk.* Choose $k < n \in \mathbb{Z}$. Define

$$(2) \qquad \begin{aligned} C_k^n = \{&w = (w(k), \ldots, w(n)); w(l) \in \mathbb{Z}^d \text{ for } l \in \{k, \ldots, n\} \\ &\text{and } |w(l+1) - w(l)|^2 = d \text{ for } l \in \{k, \ldots n-1\}\}, \end{aligned}$$

where $|\cdot|$ is the $L^2$-norm in $\mathbb{R}^d$. For $x, y \in \mathbb{Z}^d$, we define

(3) $$C_{k,x}^n = \{w \in C_k^n;\ w(k) = x\},$$

(4) $$C_k^{n,y} = \{w \in C_k^n;\ w(n) = y\}.$$

Then $P_{k,x}^n$ and $P_k^{n,y}$ denote the uniform probability measures on the sets $C_{k,x}^n$ and $C_k^{n,y}$, respectively. $P_{k,x}^n$ is the $d$-dimensional random walk starting at time $k$ in site $x$; $P_k^{n,y}$ is the random walk measure in reverted time. The expectations are denoted by $E_{k,x}^n$ and $E_k^{n,y}$, respectively. Hence, for $w \in C_{k,x}^n$,

(5) $$(S_i)_{i=k,\ldots,n} = (i, w(i))_{i=k,\ldots,n}$$

is a directed random walk under $P_{k,x}^n$ in $\mathbb{Z}^{d+1}$.

3. *Random monomers.* $\omega = (\omega_i)_{i \in \mathbb{Z}}$: an i.i.d. sequence of random variables taking the values $\pm 1$ with probability $1/2$ each; $\mathbb{P}_\omega$ denotes its probability law. Furthermore, we denote by $\mathbb{P}$ the product measure of $\mathbb{P}_\eta$ and $\mathbb{P}_\omega$; $\Omega$ denotes the space of all couples $(\omega, \eta)$ and $\mathbb{E}$ denotes the mean w.r.t. $\mathbb{P}$.

*Heteropolymer.* Now we can define the heteropolymer: Fix $\lambda \in [0, \infty)$, $h \in \mathbb{R}$, $x \in \mathbb{Z}^d$ and $k < n \in \mathbb{Z}$. Given $(\omega, \eta) \in \Omega$, define the transformed probability law on the path space $C_{k,x}^n$ by putting

(6) $$Q_{k,x}^n(w)(\omega, \eta) = \frac{1}{Z_{k,x}^n(\omega, \eta)} \exp\left\{\lambda \sum_{i=k+1}^n \Delta_\eta(S_i)(\omega_i + h)\right\} P_{k,x}^n(w),$$

where $S_i = (i, w(i))$, $Z_{k,x}^n(\omega, \eta)$ is the normalizing partition sum and

(7) $$\Delta_\eta(S_i) = \begin{cases} -1, & \text{if } S_i = (i, w(i)) \in D_\eta, \\ +1, & \text{otherwise.} \end{cases}$$

Analogously, we define $Q_k^{n,x}(w)(\omega, \eta)$. By the definition of our model, we have that $Z_{k,x}^n \stackrel{(d)}{=} Z_k^{n,x}$.

*Interpretation.* We think of $(S_i)_{i \geq k} = (i, w(i))_{i \geq k}$ as a directed polymer in $\mathbb{Z}^{d+1}$, starting at $(k, x)$, consisting of monomers represented by the sites in the path. This polymer is in a random medium described by $D_\eta$. One should think that one has, for example, water with random oil droplets $D_\eta$. The monomers are now of two different types, occurring in a random order indexed by $\omega$. Namely, $\omega_i = +1$ means that monomer $i$ is hydrophilic, $\omega_i = -1$ that it is hydrophobic. Since $\Delta_\eta(S_i) = +1$ when monomer $i$ lies in the water and $\Delta_\eta(S_i) = -1$ when it lies in the oil, we see that the weight factor in (6) encourages matches and discourages mismatches for the first $n$ monomers. For $h = 0$, both types of monomers interact equally strongly with the water and with the oil. For $h \neq 0$, on the other hand, the interaction strength is asymmetric. The parameter $\lambda$ is the overall interaction strength and plays the role of inverse temperature.

1.2. *Main results.*  First we prove the existence of the infinite volume limit of the free energy:

THEOREM 1.1.  *For every $\lambda \geq 0$, $h \in \mathbb{R}$ and $p \in [0,1]$, the following limit exists and*

$$\text{(8)} \qquad \lim_{n \to \infty} \frac{1}{n} \log Z_{0,0}^n = \Phi_p(\lambda, h), \qquad \mathbb{P}\text{-}a.s. \text{ and } L^1(\mathbb{P}),$$

*where $\Phi_p(\lambda, h) \geq \lambda h$ is a deterministic number.*

One main property of our model is to obtain a phase transition picture. We prove that there exists a critical curve which divides the parameter space into to regimes, namely, the localized and the delocalized phase.

DEFINITION 1.2.  We say the polymer is:

(a) *localized* if $\Phi_p(\lambda, h) > \lambda h$,
(b) *delocalized* if $\Phi_p(\lambda, h) = \lambda h$.

THEOREM 1.3.  *For every $p \in (0,1]$ and $\lambda > 0$, there exists $h_p^{(c)}(\lambda) \in \mathbb{R}$ such that the polymer is:*

(a) *localized for $h < h_p^{(c)}(\lambda)$,*
(b) *delocalized for $h \geq h_p^{(c)}(\lambda)$.*

*Moreover, $\lambda \mapsto h_p^{(c)}(\lambda)$ is continuous and nondecreasing on $(0, \infty)$.*

The nice feature of our model is that we are also able to give polymer path properties in both regimes.

We define for $d \geq 3$ (random walk is transient for $d \geq 3$):

$$\text{(9)} \qquad \alpha(d) = P_{0,0}^{\infty}[w(i) \text{ does not reenter state } 0 \text{ for all } i \geq 1] > 0.$$

THEOREM 1.4 (Delocalized behavior).  *Choose $d \geq 3$, $\lambda > 0$ and $h \geq (\frac{1}{2\lambda} \log \cosh(2\lambda)) \vee (1 - \frac{1}{2\lambda} \log \frac{1}{1-\alpha(d)})$. Then we have, for all $a_0 > 0$,*

$$\text{(10)} \qquad \liminf_{n \to \infty} Q_{0,0}^n[|w(n)| > a_0 n^{1/2}] > 0, \qquad \mathbb{P}\text{-}a.s.,$$

$$\text{(11)} \qquad \limsup_{c_0 \to \infty} \limsup_{n \to \infty} Q_{0,0}^n[|w(n)| > c_0 n^{1/2}] = 0, \qquad \mathbb{P}\text{-}a.s.$$

REMARK.  Of course, the choice $(\lambda, h)$ belongs to the delocalized regime in the above theorem [see Lemma 2.3 and (31)]. Hence, in Theorem 1.4 we prove that the polymer behaves diffusively in the large $h$ regime for $d \geq 3$, that is, the directed polymer lives on the same scale as the directed random walk, namely, on $n^{1/2}$.

Next we state that the endpoint of the polymer is exponentially tight in the localized regime.

THEOREM 1.5 (Localized regime). *Choose $d \geq 1$ and $(\lambda, h) \in \{(\lambda, h) \in (0, \infty) \times \mathbb{R}; \Phi_p(\lambda, h) > \lambda h\}$: There is a $\mathbb{P}$-a.s. finite random variable $\nu : \Omega \to \mathbb{N}$, $\varepsilon > 0$ and $c > 0$ such that $\mathbb{P}$-a.s. for all large $n$,*

$$(12) \qquad Q^0_{-n,0}[w(0) = z] \leq c \exp\{-\varepsilon |z|\}, \qquad |z| > \nu(\omega, \eta).$$

Define on $C^n_{0,0}$ the following annealed measure:

$$(13) \qquad \mathbb{Q}_n(dw) = \int \mathbb{P}(d\omega, d\eta) Q^n_{0,0}(dw)(\omega, \eta).$$

COROLLARY 1.6. *The laws of $w(n)$ under $\mathbb{Q}_n$, $n \geq 0$, are tight.*

1.3. *Interpretation and one-interface heteropolymer.* The model defined in (6) is one step into the direction of considering a heteropolymer in a random medium. The random medium is chosen such that the droplets are on the first coordinate axis. For $d = 1$, this has a good interpretation: The model is equivalent to the model where the directed random walk is reflected at the $x$-axis. This way one can interpret the medium such that we have water with random oil droplets at the surface.

More general geometries for the random droplet can be obtained in a similar way as in [6] (folding procedure). Then the droplets are located in parallel layers, where the distances of the layers are such that the heteropolymer does not see too many droplets at the same time.

Since there is an increasing interest into these models, we conclude the first section with some remarks on the related models. A one-interface model is studied by Garel, Huse, Leibler and Orland [4], Albeverio and Zhou [1] and Bolthausen and den Hollander [3], and early studies include Sinai [12] ($h = 0$) and Grosberg, Izrailev and Nechaev [5] ($\omega$ periodic). Recent results on related one-interface models appear in [9, 10, 11].

These authors have studied related models: The medium has two different solvents which are separated by a linear interface. They prove that in this situation there is a phase transition depending on the coefficients $h$ and $\lambda$ (we see that this is also true in our model). The two regimes are called localized regime and delocalized regime.

In most of these models it is much more difficult to prove results on path behaviors, than on free energies. For example, in [2] various path properties are derived in the localized regime, but there is not much proved in the delocalized regime, for example, there is no proof for the expected diffusive path behavior. In our model we have both results, tightness in the localized regime and, moreover, we prove diffusive behavior for some part of the delocalized regime for $d \geq 3$. The reason for the fact that we can only prove diffusive behavior for $d \geq 3$ is the same as in [13], namely, that the usual random walk is transient for high dimensions. We do not know how to control the occasional returns of the random walk to the origin in $d = 1, 2$.

A first step into the direction of relaxing the one-interface model is done in [6]. They have studied the path behavior of a heteropolymer in a medium with two different solvents arranged in alternating layers. Another model which comes away from the linear interface is studied in [7]. They study a medium where the two solvents are arranged in large blocks. In most of these models one cannot prove much about path properties. For explicit open challenges, we refer to Section 1.5 in [6] and Section 1.7 in [7].

**2. Phase transition picture.** We define

$$(14) \quad \widehat{Z}_{k,x}^n(\omega, \eta) = E_{k,x}^n \left[ \exp\left\{ \lambda \sum_{i=k+1}^n \Delta_\eta(S_i)(\omega_i + h) \right\}, w(n) = 0 \right];$$

this is the partition sum $Z_{k,x}^n$ restricted to end at the origin at time $n$, that is, $w(n) = 0$.

LEMMA 2.1. *The following limit exists and*

$$(15) \quad \lim_{n \to \infty} \frac{1}{2n} \log \widehat{Z}_{0,0}^{2n} = \Phi_p(\lambda, h), \qquad \mathbb{P}\text{-}a.s. \text{ and } L^1(\mathbb{P}),$$

*where $\Phi_p(\lambda, h)$ is a deterministic number.*

PROOF. The proof is similar to the proof of Lemma 1 in [3]. It uses Kingman's [8] subadditive ergodic theorem in its standard way. □

PROOF OF THEOREM 1.1. The lower bound easily follows from the fact that $Z_{0,0}^{2n} \geq \widehat{Z}_{0,0}^{2n}$, hence, it remains to prove the upper bound. We define the following events:

$$(16) \qquad A_{k,n} = \{w \in C_{2k}^{2n}; \ w(i) \neq 0 \text{ for } 2k < i \leq 2n\},$$

$$(17) \qquad B_{k,n} = \{w \in C_{2k}^{2n}; \ w(i) \neq 0 \text{ for } 2k < i < 2n, w(2n) = 0\}.$$

Furthermore, $a_{n-k} = P_{2k,0}^{2n}[A_{k,n}]$ and $b_{n-k} = P_{2k,0}^{2n}[B_{k,n}]$. Now we have

$$Z_{0,0}^{2n} = \widehat{Z}_{0,0}^{2n} + \sum_{k=0}^{n-1} \widehat{Z}_{0,0}^{2k} \cdot E_{2k,0}^{2n}\left[\exp\left\{\lambda \sum_{i=2k+1}^{2n} \Delta_\eta(S_i)(\omega_i + h)\right\}, A_{k,n}\right]$$

$$(18) \qquad \leq \widehat{Z}_{0,0}^{2n} + \sum_{k=0}^{n-1} \widehat{Z}_{0,0}^{2k} \cdot E_{2k,0}^{2n}\left[\exp\left\{\lambda \sum_{i=2k+1}^{2n} \Delta_\eta(S_i)(\omega_i + h)\right\}, B_{k,n}\right]$$

$$\times \frac{a_{n-k}}{b_{n-k}} \cdot \exp\{\lambda(1 + |h|)\}.$$

As in [3], Lemma 2, the claim of the theorem follows from the fact that there exists $c_1 > 0$ such that, for all $k \in \{1, \ldots, n\}$,

$$\frac{a_k}{b_k} \leq c_1 n^d, \tag{19}$$

hence, we obtain

$$Z_{0,0}^{2n} \leq (1 + c_1 n^d \exp\{\lambda(1 + |h|)\}) \widehat{Z}_{0,0}^{2n}, \tag{20}$$

which finishes, in view of Lemma 2.1, the infinite volume statement. It remains to prove the lower bound on $\Phi_p$: Let $w^{(k)}(i)$ denote the $k$th coordinate of $w(i)$:

$$\frac{1}{n} \log Z_{0,0}^n$$

$$\geq \frac{1}{n} \log E_{0,0}^n \left[ \exp\left\{ \lambda \sum_{i=1}^n \Delta_\eta(S_i)(\omega_i + h) \right\}, \right.$$

$$\left. w^{(k)}(i) \neq 0, i = 1, \ldots, n, k = 1, \ldots, d \right] \tag{21}$$

$$= \frac{\lambda}{n} \sum_{i=1}^n (\omega_i + h)$$

$$+ \frac{1}{n} \log P_{0,0}^n[w^{(k)}(i) \neq 0, i = 1, \ldots, n, k = 1, \ldots, d].$$

The first term on the right-hand side of the above inequality converges $\mathbb{P}$-a.s. to $\lambda h$, whereas the last term converges to 0 because

$$P_{0,0}^n[w^{(k)}(i) \neq 0, i = 1, \ldots, n, k = 1, \ldots, d] = O(n^{-3d/2}). \tag{22}$$

This completes the proof of Theorem 1.1. □

DEFINITION 2.2. We define the following sets:

$$\mathcal{L}_p = \{(\lambda, h) \in (0, \infty) \times \mathbb{R}; \ \Phi_p(\lambda, h) > \lambda h\}, \tag{23}$$

$$\mathcal{D}_p = ((0, \infty) \times \mathbb{R}) \setminus \mathcal{L}_p. \tag{24}$$

LEMMA 2.3. *For $p \in (0, 1]$, both $\mathcal{L}_p$ and $\mathcal{D}_p$ are nonempty.*

PROOF. *$\mathcal{L}_p$ is nonempty.* We mention that

$$\Phi_p(\lambda, h) > \lambda h$$

$$\iff \quad \mathbb{P}\text{-a.s. and } L^1(\mathbb{P}): \tag{25}$$

$$\lim_{n \to \infty} \frac{1}{n} \log E_{0,0}^n \left[ \exp\left\{ \lambda \sum_{i=k+1}^n (\Delta_\eta(S_i) - 1)(\omega_i + h) \right\} \right] > 0;$$

this follows from the fact that $\lim_n \frac{1}{n}\sum_{i=1}^n \omega_i = 0$, $\mathbb{P}$-a.s. and $L^1(\mathbb{P})$. Therefore, we prove that the right-hand side of (25) is true for certain choices of $\lambda$ and $h$. The main idea in the following proof is to introduce an indicator function in the expression on the right-hand side of (25). This indicator allows for a renewal argument which will give a positive lower bound.

We choose $n \geq 1$ an odd number and assume that $1 = (1,\ldots,1) \in \mathbb{Z}^d$ is the one-vector. Then we have, using translation invariance in the second step,

$$\mathbb{E}\left[\frac{1}{n}\log E_{0,0}^n\left[\exp\left\{\lambda\sum_{i=1}^n(\Delta_\eta(S_i)-1)(\omega_i+h)\right\}\right]\right]$$

$$\geq \mathbb{E}\left[\frac{1}{n}\log E_{0,0}^n\left[\exp\left\{\lambda\sum_{i=2}^n(\Delta_\eta(S_i)-1)(\omega_i+h)\right\},w(1)=1\right]\right]$$

(26)
$$= \mathbb{E}\left[\frac{1}{n}\log E_{0,1}^{n-1}\left[\exp\left\{\lambda\sum_{i=1}^{n-1}(\Delta_\eta(S_i)-1)(\omega_i+h)\right\}\right]\right]$$

$$\geq \mathbb{E}\left[\frac{1}{n}\log E_{0,1}^{n-1}\left[\exp\left\{\lambda\sum_{i=1}^{n-1}(\Delta_\eta(S_i)-1)(\omega_i+h)\right\},\right.\right.$$

$$\left.\left.\bigcap_{i=1}^{(n-1)/2}\{w(2i)=1\}\right]\right].$$

Using the Markov property and translation invariance, we can decouple the last term and obtain

(27)
$$\mathbb{E}\left[\frac{1}{n}\log E_{0,0}^n\left[\exp\left\{\lambda\sum_{i=1}^n(\Delta_\eta(S_i)-1)(\omega_i+h)\right\}\right]\right]$$

$$\geq \frac{1}{n}\sum_{i=1}^{(n-1)/2}\mathbb{E}[\log E_{0,1}^2[\exp\{\lambda(\Delta_\eta(S_1)-1)(\omega_1+h)\},w(2)=1]].$$

An easy lower bound for this last term can be obtained as follows: Assume $2 = (2,\ldots,2) \in \mathbb{Z}^d$ is the "two-vector," then

$$\mathbb{E}[\log E_{0,1}^2[\exp\{\lambda(\Delta_\eta(S_1)-1)(\omega_1+h)\},w(2)=1]]$$

$$\geq \mathbb{E}[\log E_{0,1}^2[\exp\{\lambda(\Delta_\eta(S_1)-1)(\omega_1+h)\},$$

$$w(1)\in\{0,2\},w(2)=1]]$$

(28)
$$= \mathbb{E}\left[\log\left(\frac{1}{4^d}+\frac{1}{4^d}\exp\{-2\lambda(\omega_1+h)\mathbf{1}_{\{\eta_1=1\}}\}\right)\right]$$

$$= -d\log 4 + (1-p)\log 2 + p\mathbb{E}[\log(1+\exp\{-2\lambda(\omega_1+h)\})]$$

$$\geq -(2d-(1-p))\log 2 + \frac{p}{2}\log(1+\exp\{2\lambda(1-h)\})$$

$$\geq -(2d - (1-p))\log 2 + p\lambda(1-h).$$

This last term is strictly positive for

(29) $$h < 1 - \frac{(2d - (1-p))\log 2}{p\lambda}.$$

This proves that $\mathcal{L}_p$ is nonempty for $p > 0$.

$\mathcal{D}_p$ *is nonempty.* This step is similar to Step 3 in the proof of Theorem 2 of [3]. We have, using Jensen's inequality,

$$\Phi_p(\lambda, h) - \lambda h$$

$$= \lim_{n\to\infty} \mathbb{E}\left[\frac{1}{n}\log E_{0,0}^n\left[\exp\left\{\lambda \sum_{i=1}^n (\Delta_\eta(S_i) - 1)(\omega_i + h)\right\}\right]\right]$$

$$\leq \limsup_{n\to\infty} \frac{1}{n} \log E_{0,0}^n\left[\prod_{i=1}^n \mathbf{1}_{\{w(i)=0\}}\right.$$

(30) $$\left. \times \mathbb{E}[\exp\{-2\lambda(\omega_1 + h)\mathbf{1}_{\{\eta_1=1\}}\}]\right]$$

$$= \limsup_{n\to\infty} \frac{1}{n} \log E_{0,0}^n\left(\prod_{i=1}^n \mathbf{1}_{\{w(i)=0\}}\left((1-p)\right.\right.$$

$$\left.\left. + \frac{p}{2}[e^{-2\lambda(1+h)} + e^{2\lambda(1-h)}]\right)\right).$$

If the last term in the square brackets is less or equal to 2, the left-hand side is less or equal to 0. But this last condition is equivalent to

(31) $$h \geq \frac{1}{2\lambda} \log \cosh(2\lambda).$$

This proves that $\mathcal{D}_p$ is nonempty, and hence, completes the proof of Lemma 2.3. □

PROOF OF THEOREM 1.3. From Lemma 2.3, we already know that both $\mathcal{L}_p$ and $\mathcal{D}_p$ are nonempty for $p > 0$. Hence, it remains to prove the existence of a continuous and nondecreasing phase transition curve $\lambda \mapsto h_p^{(c)}(\lambda)$.

*Step* 1. From (31), we know that we are in the delocalized regime for $h \geq 1$. So assume $h < 1$. If $(\lambda, h) \in \mathcal{D}_p$, then $(\lambda + \delta, h + \varepsilon) \in \mathcal{D}_p$ for all $\delta, \varepsilon \geq 0$ with $\varepsilon \geq \delta(1-h)/\lambda$. The proof is similar to Step 1 in the proof of Theorem 1 in [3]. This proves that there exists a function $\lambda \mapsto h_p^{(c)}(\lambda)$ such that $(\lambda, h) \in \mathcal{D}_p$ for all $h \geq h_p^{(c)}(\lambda)$: we define, for $\lambda \geq 0$ and $p > 0$,

$$h_p^{(c)}(\lambda) = \inf\{h \in \mathbb{R};\ (\lambda, h) \in \mathcal{D}_p\}. \tag{32}$$

By the continuity of $\Phi_p(\lambda, h) - \lambda h$, we have $(\lambda, h_p^{(c)}(\lambda)) \in \mathcal{D}_p$.

*Step* 2. Similar as in (25), we see that $\Psi_p(\lambda, h) = \Phi_p(\lambda, h) - \lambda h \geq 0$ is a convex curve in $\lambda$ (use Hölder's inequality) with boundary value $\Psi_p(0, h) = 0$. Hence, if $(\lambda, h) \notin \mathcal{D}_p$ (for fixed $h$), then $\Psi_p(\lambda, h) > 0$. But using the convexity and positivity of $\Psi_p(h, \cdot)$, we know that $\Psi_p(\lambda + \delta, h) > 0$ for all $\delta \geq 0$. Hence, $\lambda \mapsto h_p^{(c)}(\lambda)$ is nondecreasing. Step 1 shows that the slope at $\lambda > 0$ is bounded from above by $(1 - h_p^{(c)}(\lambda))/\lambda$ which is finite. Hence, we get the continuity on $(0, \infty)$. This completes the proof of Theorem 1.3. $\square$

**3. Path behavior in the large-$h$ regime for $d \geq 3$.** In this section we prove that the polymer behaves diffusively for large $h$ when $d \geq 3$. Choose $d \geq 1$, $\lambda \in [0, \infty)$, $h \in \mathbb{R}$, $x \in \mathbb{Z}^d$, $k < n \in \mathbb{Z}$ and $(\omega, \eta) \in \Omega$. Define the random variable

$$\begin{aligned}
\Psi_{k,x}^n(\omega, \eta) &= \frac{Z_{k,x}^n(\omega, \eta)}{\exp\{\lambda \sum_{i=k+1}^n (\omega_i + h)\}} \\
&= E_{k,x}^n\left[\exp\left\{\lambda \sum_{i=k+1}^n (\Delta_\eta(S_i) - 1)(\omega_i + h)\right\}\right].
\end{aligned} \tag{33}$$

LEMMA 3.1. *Choose $h \geq \frac{1}{2\lambda} \log \cosh(2\lambda)$. For $x \in \mathbb{Z}^d$, $k \in \mathbb{Z}$, there exists a random variable $\Psi_{k,x} \in L^1(\mathbb{P})$ such that $\Psi_{k,x}^n$ converges $\mathbb{P}$-a.s. to $\Psi_{k,x}$ as $n \to \infty$. For $d \geq 3$, we have $\Psi_{k,x} > 0$, $\mathbb{P}$-a.s., and $\Psi_{k,0} \geq \alpha(d)$.*

PROOF. Choose $n > k$ and denote by $\mathcal{F}_{k,n}$ the $\sigma$-field generated by $\omega_{k+1}, \ldots, \omega_n, \eta_{k+1}, \ldots, \eta_n$. Then $(\mathcal{F}_{k,n})_{n \geq k}$ is a filtration of $\sigma$-fields and $\Psi_{k,x}^n$ is $\mathcal{F}_{k,n}$-measurable.

Using the Markov property, we obtain

$$\mathbb{E}[\Psi_{k,x}^{n+1} | \mathcal{F}_{k,n}]$$

$$= \mathbb{E}\left[\sum_{y \in \mathbb{Z}^d} E_{k,x}^n\left[\exp\left\{\lambda \sum_{i=k+1}^n (\Delta_\eta(S_i) - 1)(\omega_i + h)\right\}, S_n = (n, y)\right]\right.$$

$$\left. \times \Psi_{n,y}^{n+1}(\omega, \eta) \Big| \mathcal{F}_{k,n}\right] \tag{34}$$

$$= \sum_{y \in \mathbb{Z}^d} E_{k,x}^n\left[\exp\left\{\lambda \sum_{i=k+1}^n (\Delta_\eta(S_i) - 1)(\omega_i + h)\right\}, S_n = (n, y)\right]$$

$$\times \mathbb{E}[\Psi_{n,y}^{n+1}(\omega, \eta) | \mathcal{F}_{k,n}].$$

Next we consider

$$
\begin{aligned}
\mathbb{E}[\Psi_{n,y}^{n+1}(\omega,\eta)|\mathcal{F}_{k,n}] \\
&= \mathbb{E}[\Psi_{n,y}^{n+1}(\omega,\eta)] \\
&= P_{y,n}^{n+1}[w(n+1) \neq 0] \\
&\quad + P_{y,n}^{n+1}[w(n+1) = 0] \\
&\quad \times \mathbb{E}[\mathbf{1}_{\{\eta_{n+1}=-1\}} + \mathbf{1}_{\{\eta_{n+1}=+1\}}\exp\{-2\lambda(\omega_{n+1}+h)\}] \\
&= 1 - P_{y,n}^{n+1}[w(n+1) = 0] \\
&\quad \times p(1 - \tfrac{1}{2}(\exp\{-2\lambda(1+h)\} + \exp\{2\lambda(1-h)\})) \\
&\leq 1, \qquad \text{as soon as } \exp\{-2\lambda(1+h)\} + \exp\{2\lambda(1-h)\} \leq 2.
\end{aligned}
\tag{35}
$$

But this last condition is equivalent to the assumption on $h$ in the lemma. Hence, we obtain

$$
\begin{aligned}
\mathbb{E}[\Psi_{k,x}^{n+1}|\mathcal{F}_{k,n}] \\
&\leq \sum_{y \in \mathbb{Z}^d} E_{k,x}^n\left[\exp\left\{\lambda \sum_{i=k+1}^{n} (\Delta_\eta(S_i)-1)(\omega_i+h)\right\},\right. \\
&\qquad\qquad\qquad\qquad\qquad\qquad\left. S_n = (n,y)\right] = \Psi_{k,x}^n.
\end{aligned}
\tag{36}
$$

This proves that $\Psi_{k,x}^n$ is a supermartingale with respect to $\mathcal{F}_{k,n}$. But then the first claim of our lemma is implied by the martingale convergence theorem.

Choose $k < n < m$. Using the Markov property, we have

$$
\begin{aligned}
\Psi_{k,x}^m &= \sum_{y \in \mathbb{Z}^d} E_{k,x}^n\left[\exp\left\{\lambda \sum_{i=k+1}^{n} (\Delta_\eta(S_i)-1)(\omega_i+h)\right\},\right. \\
&\qquad\qquad\qquad\qquad\qquad\qquad\left. S_n = (n,y)\right] \cdot \Psi_{n,y}^m.
\end{aligned}
\tag{37}
$$

Letting $m \to \infty$, we find with the martingale convergence theorem (as above) that, $\mathbb{P}$-a.s.,

$$
\begin{aligned}
\Psi_{k,x} &= \sum_{y \in \mathbb{Z}^d} E_{k,x}^n\left[\exp\left\{\lambda \sum_{i=k+1}^{n} (\Delta_\eta(S_i)-1)(\omega_i+h)\right\},\right. \\
&\qquad\qquad\qquad\qquad\qquad\qquad\left. S_n = (n,y)\right] \cdot \Psi_{n,y}.
\end{aligned}
\tag{38}
$$

Assume that $\Psi_{k,x} = 0$; this implies [using (38)] that $\Psi_{k+2l,x} = 0$, $\mathbb{P}$-a.s., for all $l \geq 1$. Hence,

$$
\begin{aligned}
A_k &= \{(\omega, \eta) \in \Omega;\ \Psi_{k,x}(\omega, \eta) = 0\} \\
&\subset A_{k+2l} = \{(\omega, \eta) \in \Omega;\ \Psi_{k+2l,x}(\omega, \eta) = 0\}, \qquad \mathbb{P}\text{-a.s.}
\end{aligned}
\tag{39}
$$

Using translation invariance, we have $\mathbb{P}[A_k] = \mathbb{P}[A_{k+2l}]$, which implies that $A_k = A_{k+2l}$, $\mathbb{P}$-a.s. Moreover, $A_{k+2l}$ is $\mathcal{F}_{k+2l,\infty}$ measurable, and $A_k$ is $\bigcap_{l \geq 0} \mathcal{F}_{k+2l,\infty}$ measurable. Using Kolmogorov's 0–1-law, it follows that $\mathbb{P}[A_k] = 0$ or 1, hence, $\mathbb{P}$-a.s., either $\Psi_{k,x} = 0$ or $\Psi_{k,x} > 0$.

Assume that $d \geq 3$, then the random walk is transient. Using definition (9), we have $\alpha(d) > 0$, and henceforth,

$$
\begin{aligned}
1 &\geq Q_{0,0}^n[w(i) \neq 0 \text{ for all } i \geq 1] \\
&= \frac{\exp\{\lambda \sum_{i=1}^n (\omega_i + h)\} P_{0,0}^n[w(i) \text{ does not reenter state } 0 \text{ for } i = 1, \ldots, n]}{Z_{0,0}^n(\omega, \eta)} \\
&\geq \frac{\alpha(d)}{\Psi_{0,0}^n},
\end{aligned}
\tag{40}
$$

which implies that $\Psi_{0,0}^n \geq \alpha(d) > 0$ for all $n \geq 1$. This implies the claim for $\Psi_{0,0}$, whereas, for general starting point $(k, x)$, the proof goes analogously. This completes the proof of Lemma 3.1. $\square$

Lemma 3.1 and (40) lead to the following corollary

COROLLARY 3.2. *Assume $d \geq 3$ and $h \geq \frac{1}{2\lambda} \log \cosh(2\lambda)$. Then there exists a random variable $\Psi_{0,0} \geq \alpha(d)$ with $\mathbb{E}[\Psi_{0,0}] < \infty$ and such that*

$$
\lim_{n \to \infty} Q_{0,0}^n[w(i) \neq 0 \text{ for all } i \geq 1] = \frac{\alpha(d)}{\Psi_{0,0}} > 0, \qquad \mathbb{P}\text{-a.s.}
\tag{41}
$$

REMARK. This means that, with positive probability, we do not return to the oil droplets, which is a kind of a transience statement.

For $w \in C_{0,0}^n$, we define the number of returns to the origin

(42)    $N_n =$ number of returns of $w(i)$ to the origin within $(0, n]$.

Let $E_{Q_{0,0}^n}$ denote the expectation w.r.t. $Q_{0,0}^n$.

PROPOSITION 3.3. *Choose $d \geq 3$, $\lambda > 0$ and $h > 1 - \frac{1}{2\lambda} \log \frac{1}{1-\alpha(d)}$, then we have*

$$
\limsup_{n \to \infty} E_{Q_{0,0}^n}[N_n] < \infty.
\tag{43}
$$

PROOF. Choose $(\omega, \eta) \in \Omega$ and $k, n \geq 1$. Then we have from (40) that

$$Q_{0,0}^n[N_n = k] = \sum_{1 \leq z_1 < \cdots < z_k \leq n} Q_{0,0}^n\left[\bigcap_{i=1}^k \{w(z_i) = 0\}, N_n = k\right]$$

$$= \frac{\exp\{\lambda \sum_{i=1}^n (\omega_i + h)\}}{Z_{0,0}^n(\omega, \eta)}$$

(44)
$$\times \sum_{1 \leq z_1 < \cdots < z_k \leq n} \left(P_{0,0}^n\left[\bigcap_{i=1}^k \{w(z_i) = 0\}, N_n = k\right]\right.$$

$$\left.\times \exp\left\{-2\lambda \sum_{i=1}^k (\omega_{z_i} + h)\mathbf{1}_{\{\eta_{z_i} = 1\}}\right\}\right)$$

$$\leq \frac{P_{0,0}^n[N_n = k]}{\alpha(d)} \cdot \frac{\alpha(d)}{\Psi_{0,0}^n} \cdot \exp\{-2\lambda \min\{h-1, 0\}\}^k$$

$$\leq P_{0,0}^n[N_n = k] \cdot \frac{1}{\alpha(d)} \cdot \exp\{-2\lambda \min\{h-1, 0\}\}^k.$$

Assume $h \geq 1$, then

(45) $$Q_{0,0}^n[N_n = k] \leq P_{0,0}^n[N_n = k]\frac{1}{\alpha(d)}.$$

This implies, for all $n \geq 1$,

(46)
$$E_{Q_{0,0}^n}[N_n] = \sum_{k \geq 0} Q_{0,0}^n[N_n \geq k] \leq \frac{1}{\alpha(d)} \sum_{k \geq 0} P_{0,0}^n[N_n \geq k]$$

$$\leq \frac{1}{\alpha(d)} \sum_{k \geq 0} P_{0,0}^\infty[N_\infty \geq k] = \frac{1}{\alpha(d)} E_{0,0}^\infty[N_\infty] < \infty.$$

This finishes the proof for $h \geq 1$.

So assume $1 > h > 1 - \frac{1}{2\lambda}\log\frac{1}{1-\alpha(d)}$ and define $\gamma = \exp\{2\lambda(1-h)\} \in (1, \frac{1}{1-\alpha(d)})$. Furthermore, we choose $\gamma_1 \in (\gamma, \frac{1}{1-\alpha(d)})$ (hence, $\gamma_1 > 1$). There exists $K_0$ such that $k\gamma^k < \gamma_1^k$ for all $k \geq K_0$. Hence, we have, for all $n \geq 1$,

$$E_{Q_{0,0}^n}[N_n] = \sum_{k \geq 0} k Q_{0,0}^n[N_n = k]$$

$$\leq K_0^2 + \frac{1}{\alpha(d)} \sum_{k \geq K_0} k\gamma^k P_{0,0}^n[N_n = k]$$

(47) $$\leq K_0^2 + \frac{1}{\alpha(d)} \sum_{k \geq K_0} \gamma_1^k P_{0,0}^\infty[N_n = k]$$

$$\leq K_0^2 + \frac{1}{\alpha(d)} E_{0,0}^\infty [\gamma_1^{N_n}]$$

$$\leq K_0^2 + \frac{1}{\alpha(d)} E_{0,0}^\infty [\gamma_1^{N_\infty}]$$

$$= K_0^2 + \frac{1}{\alpha(d)} \sum_{k \geq 0} \gamma_1^k (1 - \alpha(d))^k \alpha(d)$$

$$= K_0^2 + \frac{1}{1 - \gamma_1 (1 - \alpha(d))} < \infty.$$

This completes the proof of Proposition 3.3. □

PROOF OF THEOREM 1.4. Choose $h \geq (\frac{1}{2\lambda} \log \cosh(2\lambda)) \vee (1 - \frac{1}{2\lambda} \times \log \frac{1}{1-\alpha(d)})$.

*Proof of the diffusive lower bound.* Choose $a_0 > 0$, then we have from the central limit theorem

(48) $$\liminf_{n \to \infty} P_{0,0}^\infty [|w(n)| > a_0 n^{1/2}] = 1 - \chi_d^2(a_0^2) = a_1 > 0,$$

where $\chi_d^2$ denotes the quantiles of the $\chi^2$-distribution with $d$ degrees of freedom. Since the random walk is transient for $d \geq 3$, we can choose $K_0 < \infty$ and $a_2 \in (0, a_1)$ such that

(49) $$P_{0,0}^\infty [N_\infty \leq K_0] \geq 1 - a_2 > 1 - a_1,$$

from which we obtain

(50) $$\liminf_{n \to \infty} P_{0,0}^\infty [|w(n)| > a_0 n^{1/2}, N_\infty \leq K_0] = a_3 \geq a_1 - a_2 > 0.$$

Hence, we obtain

$$\liminf_{n \to \infty} Q_{0,0}^n [|w(n)| > a_0 n^{1/2}]$$

$$\geq \liminf_{n \to \infty} Q_{0,0}^n [|w(n)| > a_0 n^{1/2}, N_n \leq K_0]$$

(51) $$\geq \liminf_{n \to \infty} \frac{1}{\Psi_{0,0}^n} P_{0,0}^\infty [|w(n)| > a_0 n^{1/2}, N_n \leq K_0]$$

$$\times \exp\{-2\lambda K_0 (1 + |h|)\}$$

$$\geq \frac{a_3}{\Psi_{0,0}} \exp\{-2\lambda K_0 (1 + |h|)\} > 0, \qquad \mathbb{P}\text{-a.s.}$$

This finishes the proof of the lower bound.

*Proof of the diffusive upper bound.* Choose $c_0 > 0$. Then we have

$$Q_{0,0}^n[|w(n)| > c_0 n^{1/2}]$$

(52)
$$= \sum_{k \geq 0} Q_{0,0}^n[|w(n)| > c_0 n^{1/2}, N_n = k]$$

$$\leq \frac{1}{\Psi_{0,0}^n} \sum_{k \geq 0} P_{0,0}^\infty[|w(n)| > c_0 n^{1/2}, N_n = k]$$

$$\times \exp\{2\lambda k \max\{1 - h, 0\}\}.$$

If $h \geq 1$, then of course we have, $\mathbb{P}$-a.s.,

(53) $$\limsup_{c_0 \to \infty} \limsup_{n \to \infty} Q_{0,0}^n[|w(n)| > c_0 n^{1/2}] \leq \frac{1}{\Psi_{0,0}} \limsup_{c_0 \to \infty} (1 - \chi_d^2(c_0^2)) = 0,$$

which finishes the proof for $h \geq 1$. So let us assume $h \in (\frac{1}{2\lambda} \log \cosh(2\lambda) \vee 1 - \frac{1}{2\lambda} \log \frac{1}{1-\alpha(d)}, 1)$ and define $\gamma = \exp\{2\lambda(1-h)\} \in (1, \frac{1}{1-\alpha(d)})$. Furthermore, we choose $q \in (1, \frac{\log 2 - \log(1-\alpha(d))}{\log 2 + \log \gamma})$. Hence, we have, using Hölder's inequality,

$$\sum_{k \geq 0} P_{0,0}^\infty[|w(n)| > c_0 n^{1/2}, N_n = k] \exp\{2\lambda k \max\{1 - h, 0\}\}$$

(54)
$$\leq P_{0,0}^\infty[|w(n)| > c_0 n^{1/2}]^{1-1/q} \sum_{k \geq 0} P_{0,0}^\infty[N_n = k]^{1/q} \gamma^k$$

$$= P_{0,0}^\infty[|w(n)| > c_0 n^{1/2}]^{1-1/q}$$

$$\times 2 \cdot \sum_{k \geq 0} (P_{0,0}^\infty[N_n = k](2\gamma)^{kq})^{1/q} (\tfrac{1}{2})^k \tfrac{1}{2}.$$

The last term in the equation above describes a geometric distribution: Assume $X$ is geometrically distributed with $p = 1/2$ and define $Y = P_{0,0}^\infty[N_n = X](2\gamma)^{Xq}$. Hence, we have, using Jensen's inequality,

$$\sum_{k \geq 0} (P_{0,0}^\infty[N_n = k](2\gamma)^{kq})^{1/q} (\tfrac{1}{2})^k \tfrac{1}{2}$$

(55)
$$= E[Y^{1/q}] \leq E[Y]^{1/q}$$

$$= \left( \sum_{k \geq 0} P_{0,0}^\infty[N_n = k](2\gamma)^{kq} (\tfrac{1}{2})^k \tfrac{1}{2} \right)^{1/q}.$$

So we calculate this last expectation. Define $\gamma_2 = \gamma^q 2^{q-1}$. By the choice of

$q$, we have $\gamma_2 \in (\gamma, \frac{1}{1-\alpha(d)})$. Hence, it can be bounded by

$$\sum_{k\geq 0} P_{0,0}^\infty[N_n = k]\gamma_2^k \frac{1}{2} = \frac{1}{2} E_{0,0}^\infty[\gamma_2^{N_n}]$$

(56)
$$\leq \frac{1}{2} E_{0,0}^\infty[\gamma_2^{N_\infty}]$$

$$= \frac{1}{2} \frac{\alpha(d)}{1 - \gamma_2(1 - \alpha(d))} < \infty.$$

We collect now all the pieces and we obtain, $\mathbb{P}$-a.s.,

(57)
$$\limsup_{n\to\infty} Q_{0,0}^n[|w(n)| > c_0 n^{1/2}]$$

$$\leq \frac{2^{1-1/q}}{\Psi_{0,0}} (1 - \chi_d^2(c_0^2))^{1-1/q} \left(\frac{\alpha(d)}{1 - \gamma_2(1 - \alpha(d))}\right)^{1/q}.$$

Since $\lim_{c_0 \to \infty} \chi_d^2(c_0^2) = 1$, the claim of Theorem 1.4 follows. □

**4. Localized regime.** First we prove a result which is similar to the result obtained by Biskup and den Hollander [2], Lemma 3, which is their analogon in the one-interface model.

LEMMA 4.1 (Localized regime). *Let $\lambda, \delta > 0$ and $h \in \mathbb{R}$ such that $\Phi_p(\lambda, h) > \lambda h + \delta$. Choose $\varepsilon \in (0, \delta)$, then there exists $\tilde{\delta} > 0$ such that, for all large $m$,*

(58)
$$\mathbb{P}\left[\frac{1}{m} \log Z_{-m,0}^0 < \Phi_p(\lambda, h) - \varepsilon\right] \leq \exp\{-\tilde{\delta} m\}.$$

PROOF. Choose $\varepsilon \in (0, \delta)$, then we have from Theorem 1.1 and Lemma 2.1 that, for all large $n$,

(59)
$$\frac{1}{2n}\mathbb{E}[\log Z_{-2n,0}^0] \geq \Phi_p(\lambda, h) - \varepsilon/2,$$

(60)
$$\frac{1}{2n}\mathbb{E}[\log \widehat{Z}_{-2n,0}^0] \geq \Phi_p(\lambda, h) - \varepsilon/2.$$

Choose $N$ even and $k \geq 1$. As in Lemma 2.1, we have (adding additional hitting points of the first coordinate axis)

(61)
$$Z_{-kN,0}^0 \geq \prod_{j=1}^k \widehat{Z}_{-jN,0}^{-(j-1)N}.$$

From (60), we have, for all large $N$ (using translation invariance),

(62)
$$\frac{1}{N}\mathbb{E}[\log \widehat{Z}_{-jN,0}^{-(j-1)N}] \geq \Phi_p(\lambda, h) - \varepsilon/2,$$

and we know that for fixed $N$ the random variables $\log \widehat{Z}_{-jN,0}^{-(j-1)N}$ are bounded i.i.d. random variables for the running index $j = 1, 2, \ldots$. Using a standard large deviations estimate, we find the following result: for all $\varepsilon \in (0, \delta)$, there exists $N_0$ such that, for all $N \geq N_0$, there exist $c, \delta_1 > 0$ such that, for all large $k$,

$$\mathbb{P}\left[\frac{1}{kN} \log Z_{-kN,0}^0 < \Phi_p(\lambda, h) - \varepsilon\right]$$

(63)
$$\leq \mathbb{P}\left[\frac{1}{k} \sum_{j=1}^{k} \frac{1}{N} \log \widehat{Z}_{-jN,0}^{-(j-1)N} < \Phi_p(\lambda, h) - \varepsilon\right]$$

$$\leq c \exp\left\{-\frac{\delta_1}{N} \cdot kN\right\}.$$

Hence, the lemma is proved for $kN$ if we choose $\tilde{\delta} \in (0, \delta_1/N)$ ($N \geq N_0$ fixed) and $k$ large. For general $m = kN + l$, $l \in \{0, 1, \ldots, N-1\}$, the claim follows similarly:

(64)
$$Z_{-m,0}^0 \geq \prod_{j=1}^{k} \widehat{Z}_{-jN,0}^{-(j-1)N} \widehat{Z}_{-m,0}^{-m+l},$$

the last term $\widehat{Z}_{-m,0}^{-m+l}$ is small w.r.t. large $k$, that is, $\frac{1}{k} \log \widehat{Z}_{-m,0}^{-m+l}$ is arbitrarily small for large $k$. This completes the proof of Lemma 4.1 □

PROOF OF COROLLARY 1.6. Choose $\nu, c, \varepsilon$ as in Theorem 1.5. We define, for $n \geq 0$,

(65)
$$A_n = \{(\omega, \eta);\ \text{for all } m \geq n : Q_{-m,0}^0[w(0) = z](\omega, \eta) \leq c \exp\{-\varepsilon|z|\}$$
$$\text{for all } z \text{ with } |z| > \nu(\omega, \eta)\}.$$

Then $A_n \to \Omega$, $\mathbb{P}$-a.s., using Theorem, 1.5. Therefore, for $M > 0$, $n \geq 0$, we have

(66)
$$\mathbb{Q}_n(|w(n)| > 2M)$$
$$= \int \mathbb{P}(d\omega, d\eta) Q_{-n,0}^0(|w(0)| > 2M)(\omega, \eta)$$
$$\leq \mathbb{P}[\nu > M] + \mathbb{P}[A_n^c] + \mathbb{E}[\nu \leq M, A_n, Q_{-n,0}^0(|w(0)| > 2M)]$$
$$\leq \mathbb{P}[\nu > M] + \mathbb{P}[A_n^c] + c \sum_{|z|>2M} \exp\{-\varepsilon|z|\}.$$

Since $\nu$ is finite $\mathbb{P}$-a.s., $\lim_{n \to \infty} \mathbb{P}[A_n^c] = 0$ and since the last term in (66) converges to 0 as $M \to \infty$, we see that for all $\varepsilon > 0$ there exists $M > 0$ such that, for all $n$,

(67)
$$\mathbb{Q}_n(|w(n)| > 2M) \leq \varepsilon,$$

which proves tightness of $w(n)$ under $\mathbb{Q}_n$. This completes the proof. □

PROOF OF THEOREM 1.5. We define $L = \sup\{k \leq 0;\ w(k) = 0\}$, which is the last hitting time of the first coordinate axis. In the following sums we only consider the terms which have the correct parity. Hence, we have, for $z \in \mathbb{Z}^d$, $z \neq 0$,

$$Q^0_{-n,0}[w(0) = z]$$

$$= \sum_{t=|z|_\infty}^{n} Q^0_{-n,0}[L = -t, w(0) = z]$$

$$= \sum_{t=|z|_\infty}^{n} \frac{\widehat{Z}^t_{-n,0}}{Z^0_{-n,0}} \cdot E^0_{-t,0}\left[\exp\left\{\lambda \sum_{i=-t+1}^{0} \Delta_\eta(S_i)(\omega_i + h)\right\},\right.$$

(68)
$$\left. w(0) = z, w(i) \neq 0, i = -t+1,\ldots,1\right]$$

$$\leq \sum_{t=|z|_\infty}^{n} \frac{1}{Z^0_{-t,0}} \cdot E^0_{-t,0}\left[\exp\left\{\lambda \sum_{i=-t+1}^{0} \Delta_\eta(S_i)(\omega_i + h)\right\},\right.$$

$$\left. w(0) = z, w(i) \neq 0, i = -t+1,\ldots,1\right]$$

$$\leq \sum_{t=|z|_\infty}^{n} \frac{\exp\{\lambda \sum_{i=-t+1}^{0}(\omega_i + h)\}}{Z^0_{-t,0}}.$$

From Lemma 4.1, we have, using a Borel–Cantelli argument: there exists $\Omega_1 \subset \Omega$ with $\mathbb{P}[\Omega_1] = 1$ such that, for all $(\omega, \eta) \in \Omega_1$ and all $m$ sufficiently large,

(69) $$Z^0_{-m,0} \geq \exp\{m(\Phi_p(\lambda, h) - \varepsilon)\}.$$

For the other term in (68), we have that

(70) $$\frac{1}{t}\lambda \sum_{-t+1}^{0}(\omega_i + h) \xrightarrow{t \to \infty} \lambda h, \qquad \mathbb{P}\text{-a.s.}$$

Hence, $\mathbb{P}$-a.s., for $t$ large,

(71) $$\lambda \sum_{-t+1}^{0}(\omega_i + h) \leq (\lambda h + \varepsilon)t.$$

Choose $\varepsilon \in (0, \delta/4)$ recall that $\delta = \Phi_p(\lambda, h) - \lambda h > 0$ in the localized regime, then we have, $\mathbb{P}$-a.s., for large $t$,

$$(72) \qquad \frac{\exp\{\lambda \sum_{i=-t+1}^{0} (\omega_i + h)\}}{Z_{-t,0}^0} \leq \exp\{2\varepsilon t - \delta t\} \leq \exp\{-\delta t/2\}.$$

Since the last term in the above equation is integrable, we obtain from (68), $\mathbb{P}$-a.s., for all large $|z|_\infty$ and all large $n$,

$$(73) \qquad Q_{-n,0}^0[w(0) = z] \leq c_1 \exp\{-\delta |z|_\infty / 2\}.$$

Since the supremum norm $|\cdot|_\infty$ and the $L^2$-norm $|\cdot|$ are equivalent, the proof of Theorem 1.5 follows. $\square$

DEPARTMENT OF MATHEMATICS
ETH ZÜRICH
ZENTRUM HG G 32.5
CH-8092 ZÜRICH
SWITZERLAND
E-MAIL: mario.wuethrich@math.ethz.ch